\documentstyle{amsart}

\newcommand{\lan}{\langle}
\newcommand{\ran}{\rangle}

\newcommand{\ORD}{\operatorname{ORD}}

\newsymbol\restriction 1316
\newsymbol\square 1003
\newsymbol\Vdash 130D

\newcommand{\vDash}{\models}
\begin{document}

\Large

\baselineskip=.30truein

\centerline{\bf STRICT GENERICITY}

\vskip30pt

\large
   \begin{center}

Sy D. Friedman\footnote{Research supported by NSF Contract \# 92-05530.}
\end{center}

  \begin{center}
M.I.T.\end{center}

\vskip10pt

The purpose of this note is to show that unlike for set forcing, an inner
model of a class-generic extension need not itself be a class-generic
extension. Our counterexample is of the form $L[R]$, where $R$ is a real
both generic over $L$ and constructible from $O^\#.$

\vskip10pt

\noindent
{\bf Definition} \  $\lan M,A\ran,$ $M$ transitive is a {\it ground model}
if $A\subseteq M,$ $M\vDash ZFC+A$ Replacement and $M$ is the smallest
model with this property of ordinal height $ORD(M).$ $G$ is {\it literally
generic} over $\lan M,A\ran$ if for some partial -ordering $P$ definable
over $\lan M,A\ran,$ $G$ is $P$-generic over $\lan M,A\ran$ and $\lan
M[H],A,H\ran\vDash ZFC+(A,H)$-Replacement for all $P$-generic $H.$ $S$ is
{\it generic} over $M$ if for some $A$, $S$ is definable over $\lan
M[G],A,G\ran$ for some $G$ which is literally generic over $\lan M,A\ran,$
and $S$ is {\it strictly generic} over $M$ if we also require that $G$ is
definable over $\lan M[S],A,S\ran.$

The following is a classic application of Boolean-valued forcing and can be
found in Jech [?], page ?.

\vskip10pt

\noindent
{\bf Proposition 1.} \ If $G$ is $P$-generic over $\lan M,A\ran$ where $P$
is an element of $M,S$ definable over $\lan M[G],A\ran$ then $S$ is strictly
generic over $M.$

\vskip10pt

\noindent
{\bf Proof Sketch.} \ We can assume that $P$ is a complete Boolean algebra
in $M$ and that $S\subseteq \alpha $ for some ordinal $\alpha \in M.$ Then
$H=G\cap P_0$ is $P_0$-generic over $M,$ where $P_0=$ complete subalgebra
of $P$ generated by the Boolean values of the sentences ``$\hat\beta
\in\sigma$'', where $\beta <\alpha $ and $\sigma $ is a $P$-name for $S.$
Then $H$ witnesses the strict genericity of $S.$ $\dashv$

\vskip10pt

Now we specialize to the ground model $\lan L,\phi \ran,$ under the
assumption that $O^\#$ exists.

\vskip10pt

\noindent
{\bf Theorem 2.} \ There is a real $R\in L[O^\#]$ which is generic but not
strictly generic over $L$.

Our strategy for proving Theorem 2 comes from the following observation.

\vskip10pt

\noindent
{\bf Proposition 3.} \ If $R$ is a real strictly generic over $L$ then for
some $L$-amenable $A,$ Sat$(L[R])$ is definable from $R,A,$ where
Sat denotes the Satisfaction relation.

\vskip10pt

\noindent
{\bf Proof.} \ Suppose that $A,G$ witness that $R$ is strictly generic over
$L.$ Let $G$ be $P$-generic over $\lan L,A\ran,$ $P$ definable over $\lan
L,A\ran,$ $R\in L[G],$ $G$ definable over $\lan L[R],A\ran.$ Also assume
that $\lan L[H],A,H\ran\vDash ZFC+(A,H)$-replacement for all $P$-generic
$H.$ The latter implies that the Truth and Definability Lemmas hold for
$P$-forcing, by a result of $M.$ Stanley (See Stanley [?] or
Friedman[???]). Then we have: $L[R]\vDash\varphi$ iff $\exists p\in G(p\Vdash
\varphi $ holds in $L[\sigma ])$ where $\sigma $ is a $P$-name for $R$ and
therefore Sat$(L[R])$ is definable from $R,$ Sat$\lan L,A\ran.$ As $A$ is
$L$-amenable and $O^\#$ exists, Sat$\lan L,A\ran$ is also $L$-amenable.
$\dashv$

\vskip10pt

\noindent
{\bf Remarks} \ (a) Sat$(L[R])$ could be replaced by Sat$(\lan L[R],A\ran)$
in Proposition 3, however we have no need here for this stronger
conclusion.  \ (b) \ A real violating the conclusion of Proposition 3 was
constructed in Friedman [??], however the real constructed there was not
generic over $L.$

Thus to prove Theorem 2 it will suffice to find a generic $R\in L[O^\#]$
such that for each $L$-amenable $A,$ Sat$(L[R])$ is not definable (with
parameters) over $\lan L[R],A\ran.$ First we do this not with a real $R$
but with a generic class $S,$ and afterwards indicate how to obtain $R$ by
coding $S.$

We produce $S$ using the Reverse Easton iteration $P=\lan P_\alpha |\alpha \le
\infty\ran,$ defined as follows. $P_0=$trivial forcing and for limit
$\lambda \le\infty,$ Easton support is used to define $P_\lambda $ (as a
direct limit for $\lambda $ regular, inverse limit otherwise). For singular
$\alpha ,$ $P_{\alpha +1}=P_{\alpha ^*} Q{(\alpha )}$ where $Q{(\alpha
)}$ is the trivial forcing and finally for regular $\alpha ,$ $P_{\alpha +1}=
P_\alpha *Q{(\alpha )}$ where $Q(\alpha )$ is defined as follows: let
$\lan b_\gamma |\gamma <\alpha \ran$ be the $L$-least partition of the odd
ordinals $<\alpha $ into $\alpha $-many disjoint pieces of size $\alpha $
and we take a condition in $Q{(\alpha )}$ to be $p=\lan
p(0),p(1),\dots\ran$ where for some $\alpha (p)<\alpha ,$ $p(n):\alpha
(p)\longrightarrow 2$ for each $n.$  Extension is defined by: \ $p\le q$
iff $\alpha (p)\ge\alpha (q),$ $p(n)$ extends $q(n)$ for each $n$ and
$q(n+1)(\gamma )=1$, $\delta \in b_\gamma \cap [\alpha (q),\alpha
(p))\longrightarrow p(n)(\delta )=0.$ Thus if $G$ is $Q(\alpha )$-generic
and $S_n=\bigcup\{p(n)|p\in G\}$ then $S_{n+1}(\gamma )=1$ iff $S_n(\delta
)=0$ for sufficiently large $\delta \in b_\gamma.$

Now we build a special $P$-generic  $G(\le\infty),$ definably over
$L[O^\#].$ The desired generic but not strictly generic class is
$S_0=\bigcup\{p(0)|p\in G(\infty)\}.$ We define $G(\le i_\alpha )$ by
induction on $A\in\ORD,$ where $\lan i_\alpha |\alpha \in\ORD\ran$ is the
increasing enumeration of $I\cup\{0\},$ $I=$ Silver Indiscernibles. $G(\le i_0)$ is trivial and for
limit $\lambda \le \infty,$ $G(<i_\lambda )=\bigcup\{G(<i_\alpha )|\alpha
<\lambda \},$ $G(i_\lambda )=\bigcup\{G(i_{2\alpha })|\alpha <\lambda \}$
(where $i_\infty=\infty).$

Suppose that $G(\le i_\lambda )$ is defined, $\lambda $ limit or $0,$ and
we wish to define $G(\le i_{\lambda +n})$ for $0<n<\omega .$ If $n$ is even
and $G(\le i_{\lambda +n})$ has been defined then we define $G(\le
i_{\lambda +n+1})$ as follows: \ $G(<i_{\lambda +n+1})$ is the
$L[O^\#]$-least generic extending $G(\le i_{\lambda +n}).$ To define
$G(i_{\lambda +n+1})$ first form the condition $p\in Q(i_{\lambda +n+1})$
defined by: \ $\alpha (p)=i_{\lambda +n}+1,$ $p(m)\restriction i_{\lambda
+n}=G(i_{\lambda +n})(m)=\bigcup \{q(m)|q\in G(i_{\lambda +n})\}$ for all
$m,$ $p(m)(i_{\lambda +n})=1$ iff $m>n.$ Then $G(i_{\lambda +n+1})$ is the
$L[O^\#]$-least $Q(i_{\lambda +n+1})$-generic (over $L[G(<i_{\lambda
+n+1})])$ containing the condition $p.$ If $n$ is odd and $G(\le i_{\lambda
+n})$ has been defined then we define $G(\le i_{\lambda +n+1})$ as follows:\
$G(<i_{\lambda +n+1})$ is the $L[O^\#]$-least generic extending $G(\le
i_{\lambda +n})$. To define $G(i_{\lambda +n+1}),$ first form the condition
$p\in Q(i_{\lambda +n+1})$ by: \ $\alpha (p)=i_{\lambda +n},$ $p(m)(\gamma
)=G(i_{\lambda +n})(m)(\gamma )$ for $\gamma \neq i_{\lambda +n-1}$ and
$p(m)(i_{\lambda +n-1})=0$ for all $m.$ Then $G(i_{\lambda +n+1})$ is the
$L[O^\#]$-least $Q(i_{\lambda +n+1})$-generic (over $L[G(<i_{\lambda
+n+1})])$ containing the condition $p.$ This completes the definition of
$G(\le\infty).$

Now for each $i\in I\cup\{\infty\}$ and $n\in\omega $ let
$S_n(i)=\bigcup\{p(n)|p\in G(i)\}$ and $S(i)=S_0(i),$ $S=S(\infty).$ We now
proceed to show that $S$ is not strictly-generic over $L$.

\vskip10pt

\noindent
{\bf Definition.} \ For $X\subseteq \ORD,$ $\alpha \in\ORD$ and $n\in\omega
$  we say that $\alpha $ is $X-\Sigma _n$ {\it stable} if
$\lan L_\alpha [X],X\cap \alpha \ran$ is $\Sigma _n$-elementary in $\lan
L[X],X\ran.$ $\alpha $ is $X$-{\it stable} if $\alpha $ is $X-\Sigma _n$ stable
for all $n.$

\vskip10pt

\noindent
{\bf Lemma 4.} \ For $\lambda $ limit or $0$, $n$ even, $I_{\lambda +n+1}$
is not $S$-stable.

\vskip10pt

\noindent
{\bf Proof.} \ Let $i=i_{\lambda +n}$ and $j=i_{\lambda +n+1}.$ Note that
$S_m(j)$ is defined from $S(j)$ just as $S_m(\infty)$ is defined from
$S(\infty)=S.$ But $S(j)=S\cap j$ and for $M>n,$ $S_m(j)\neq S_m(\infty)$
since $i\in S_m(j),$ $i\notin S_m(\infty).$ So $j$ is not $S$-stable.
$\dashv$

\vskip10pt

\noindent
{\bf Lemma 5.} \ For $L$-amenable $A\subseteq \ORD,$ $i_{\lambda +n+1}$ is
$(S,A)-\Sigma _n$ stable for sufficiently large limit $\lambda ,$ all
$n\in\omega.$

\vskip10pt

\noindent
{\bf Proof.} \ Let $i=i_{\lambda +n+1}$ where $\lambda $ is large enough to
guarantee that $i$ is $A$-stable. For $p\in P_{i+1}=P_i*Q(i)$ and $m\in
\omega ,$ we let $(p)_m$ be obtained from $p$ by redefining $p(i)(\bar
m)=\phi$ for $\bar m>m$ and otherwise leaving $p$ unchanged.

\vskip10pt

\noindent
{\bf Claim.} \ Suppose $\varphi$ is $\Pi_m$ relative to $S(i),$ $B$ where
$B\subseteq i,$ $B\in L.$ If $p\in P_{i+1},$ $p\Vdash\varphi$ then
$(p)_m\Vdash \varphi.$

\vskip10pt

\noindent
{\bf Proof of Claim.} \ By induction on $m\ge 1.$ For $m=1,$ if the
conclusion failed then we could choose $q\le (p)_1,$ $q(<i)\Vdash\sim\varphi$
holds of $q(0),B;$ then clearly $(q)_0\Vdash\sim\varphi,$ $(q)_0$ is
compatible with $p,$ which contradicts the hypothesis that
$p\Vdash\varphi.$ Given the result for $m,$ if the conclusion failed for
$m+1$ then we could choose $q\le (p)_{m+1},$ $q\Vdash\sim\varphi.$ Now
write  $\sim \varphi$ as $\exists x\psi ,\psi $ $\Pi_m$ and we see that
by induction we may assume that $(q)_m\Vdash\psi (\hat x)$ for some $x.$
But $(q)_m,p$ are compatible and $p\Vdash \sim\exists x\psi ,$ contradiction.
$\dashv$ (Claim.)

\vskip10pt

Now we prove the lemma. Suppose $\varphi$ is $\Pi_n$ and true of
$(S(i),A\cap i).$ Choose $p\in G(\le i), p\Vdash \varphi.$ Then by the
Claim, $(p)_n\Vdash\varphi.$ As $i$ is $A$-stable, $(p)_n\Vdash \varphi$ in
$P(\le\infty).$ By construction $(p)_n$ belongs to $G(\le\infty),$ in the
sense that $(p)_n$ $(<i)\in G(<i)\subseteq G(<\infty)$ and $(p)_n(i)\in
G(\infty).$ So $\varphi$ is true of $(S,A).$ $\dashv$

\vskip10pt

\noindent
{\bf Theorem 6.} \ $S$ is generic, but not strictly generic, over $L.$

\vskip10pt

\noindent
{\bf Proof.} \ By Proposition 3 (which also holds for classes), if $S$ were
strictly generic over $L$ then for some $L$-amenable $A$ we would have that
Sat$\lan L[S],S\ran$ would be definable over $\lan L[S],S,A\ran.$ But then
for some $n,$ all sufficiently large $(S,A)-\Sigma _n$ stables would be
$S$-stable, in contradiction to Lemmas 4,5. $\dashv$

\vskip10pt

To prove Theorem 2 we must show that an $S$ as in Theorem 6 can be coded by
a real $R$ in such a way as to preserve the properties stated in lemmas
4,5. We must first refine the above construction:

\vskip10pt

\noindent
{\bf Theorem 7.} \ Let  $\lan A(i)|i\in I\ran$ be a sequence such that
$A(i)$ is a constructible subset of $i$ for each $i\in I.$ Then there
exists $S$ obeying Lemmas 4,5 such that in addition, $A(i)$ is definable
over $\lan L_i[S],S\cap i\ran$ for $i\in \text{ Odd }
(I)=\{i_{\lambda +n}|\lambda$ limit or $0,n$ odd\}.

\vskip10pt

\noindent
{\bf Proof.} \ We use a slightly different Reverse Easton iteration: \
$Q(\alpha )$ specifies $n(\alpha )\le\omega $ and if $n(\alpha )<\omega ,$
it also specifies a constructible $A(\alpha )\subseteq\alpha ;$ then
conditions and extension are as before, except we now require that if
$n(\alpha )<\omega $ then for $p$ to extend $q,$ we must have $p(n(\alpha))
(2\beta+2)=
1$ iff $\beta \in A(\alpha),$ for $2\beta+2 \in [\alpha (q),\alpha (p)).$ 
Then if $n(\alpha )<\omega ,$ the $Q(\alpha )$-generic will code $A(\alpha )$ definably (though the
complexity of the definition increases with $n(\alpha )<\omega ).$

\vskip10pt

Now in the construction of $G(\le i_\alpha ),\alpha \le \infty$ we proceed
as before, with the following additional specifications: \ $n(i_{\lambda
+n})=n$ for odd $n$ and $n(i_{\lambda +n})=\omega $ for even $n$ ($\lambda
$ limit or $0$). And for odd $n$ we specify $A(i_{\lambda+n})$ to be the
$A(i),i=i_{\lambda +n}$ as given in the hypothesis of the Theorem.

Lemma 4 holds as before; we need a new argument for Lemma 5. Note that for
$i\in \text{ Odd}(I)$ it is no longer the case that $P(<i)\Vdash
Q(i)=Q(\infty)\cap L_i[G(<i)].$ Let $Q^*(i)$ denote $Q(\infty)\cap
L_i[G(<i)],$ i.e., the forcing $Q(i)$ where $n(i)$ has been specified as
$\omega.$ Define $(p)_m$ as before for $p\in P(\le i).$

\vskip10pt

\noindent
{\bf Claim.} \ Suppose $m\le n+1,$ $n$ is even, $i=i_{\lambda +n+1}$
$(\lambda $ limit or $0)$ and $\varphi$ is $\Pi_m$ relative to $S(i),$ $B$
with parameters, where $B\subseteq i,$ $B\in L.$ If $p\in P(\le i)$ (where
$n(i)=n+1)$ then $p\Vdash\varphi$ in $P(\le i)$ iff $(p)_m\Vdash\varphi$ in
$P^*(\le i)=P(<i)*Q^*(i)$ iff $p\Vdash\varphi$ in $P^*(\le i).$

\vskip10pt

\noindent
{\bf Proof.} \ As in the proof of the corresponding Claim in the proof of
Lemma 5.  If $m=1$ and $p\Vdash \varphi$ in $P(\le i),$ then if the
conclusion failed, we could choose $q\le(p)_1$ in $P^*(\le i),$
$q\Vdash\sim\varphi;$ then (we can assume) $(q)_0\Vdash\sim\varphi$ in
$P(\le i),$ but $(q)_0$ and $p$ are compatible. The other implications are
clear, as $P(\le i)\subseteq P^*(\le i).$ Given the result for $m\le n,$
$\varphi$ $\Pi_{m+1}$ and $p\Vdash\varphi$ in $P(\le i),$ if the conclusion
failed we could choose $q\le (p)_{m+1}$ in $P^*(\le i),q\Vdash\sim\varphi$
(indeed, $q\Vdash\sim\psi (x)$ some $x,$ where $\varphi=\forall x\psi ,\psi
\Sigma _m);$ then $q\Vdash\sim\varphi$ in $P(\le i),$
$(q)_m\Vdash\sim\varphi$ in $P^*(\le i),$ $(q)_m\Vdash\sim\varphi$ in
$P(\le i)$ by induction. But $(q)_m,p$ are compatible in $P(\le i),$ using
the fact that $m\le n$ and $q\le(p)_{m+1},$ contradiction. And again the
other implications follow, as $P(\le i)\subseteq P^*(\le i).$ $\dashv$
(Claim.)

\vskip10pt

Now the proof of Lemma 5 proceeds as before, using the new version of the
Claim. $\dashv$

\vskip10pt

The choice of $\lan A(i)|i\in I\ran$ that we have in mind comes from the
next Proposition.

\vskip10pt

\noindent
{\bf Proposition 8.} \  For each $n$ let $A_n=\{\alpha |$ For $i<j_1<\dots$
$<j_n$ in $I,\alpha <i,$ $(\alpha ,j_1\dots j_n)$ and $(i,j_1\dots j_n)$ satisfy the same formulas in $L$ with parameters $<\alpha \}.$ Then any $L$-amenaable $A$ is $\Delta
_1$-definable over $\lan L,A_n\ran$ for some $n.$

\vskip10pt

\noindent
{\bf Proof.} \ For each $i\in I, A\cap L_i$ belongs to $L$ and hence is of
the form $t(i)(\vec j_0(i),i,\vec\infty(n(i)))$ where $t(i)$ is a
$\Delta_0$-Skolem term for $L,$ $\vec j_0(i)$ is a finite sequence of
indiscernibles $<i$ and $\vec\infty$ $(n(i))$ is any sequence of
indiscernibles $>i$ of length $n(i)\in\omega.$  By Fodor's Theorem and
indiscernibility we can assume that $t(i)=t,$ $\vec j_0(i)=\vec j_0$ and
$n(i)=n$ are independent of $i.$ To see that $A$ is $\Delta _1$-definable
over $\lan L,A_{n+1}\ran$ it suffices to show that for $\vec i<\vec j$
increasing sequences from $A_{n+1}$ of length $n+1$, $\vec i$ and $\vec j$
satisfy the same formulas in $L$ with parameters $<$ min $(\vec i).$  But by
definition, for $\alpha <$ min$(\vec i)$ and $\vec i=\{i_0,\dots, i_n\},$
$\vec j=\{j_0,\dots j_n\}$ we get: \ $L\vDash\varphi(\alpha ,j_0\dots
j_n)\longleftrightarrow \varphi(\alpha ,i_0,j_1\dots
j_n)\longleftrightarrow$ $\varphi(\alpha ,i_0,i_1,j_2\dots
j_n)\longleftrightarrow \dots\longleftrightarrow \varphi(\alpha ,i_0,\dots,
i_n).$ $\dashv$

\vskip10pt

Now for $i\in I$ write $i=i_{\lambda +n},\lambda $ limit or $0$,
$n\in\omega $ and let $A(i)=A_n\cap i$. Thus by Theorem 7 there is $S$
obeying Lemmas 4,5 such that $A_n\cap i$ is definable over $\lan
L_i[S],S\cap i\ran$ for $i=i_{\lambda +n+1},n$ even.

\vskip10pt

\noindent
{\bf Proof of Theorem 2} \ First observe that as in Friedman [?], we may
build $G(\le\infty)$ to satisfy Theorem 7 for the preceding choice of $\lan
A(i)|i\in I\ran$ and in addition preserve the indiscernibility of Lim $I.$
Then by the technique of Beller-Jensen-Welch [82], Theorem 0.2 we may code
$(G(<\infty),S)$ by a real $R,$ where $S=G_0(\infty).$ The resulting $R$
obeys Lemma 4 because $S$ is definable from $R;$ to obtain  Lemma 5 for
$R$ we must modify the coding of $(G(<\infty),S)$ by $R$ in the following
way: for inaccessible $\kappa $ we require that any coding condition with
$\kappa $ in its domain reduce any dense $D\subseteq P^{<\kappa
}=\{q|\alpha (q)<\kappa \}$ strictly below $\kappa ,$ when $D$ is definable
over $\lan L_\kappa [G(<\infty),S],$ $G(<\kappa ),S\cap \kappa \ran.$ This
extra requirement does not interfere with the proofs of extendibility,
distributivity for the coding conditions (see Friedman [????]).

Now to obtain Lemma 5 for $R$ argue as follows: \ Given $L$-amenable $A,$
choose $n$ and $\lambda $ large enough so that $A$ is $\Delta _1$-definable
from $A_n$ with parameters $<i_\lambda.$ Then $i_{\lambda +n+1}$ is
$(G(<\infty),S,A)-\Sigma _n$ stable. And also $A\cap i_{\lambda +n+1}$ is
definable over $\lan L_i[G(<i),S\cap i],$ $G(<i),S\cap i\ran$ where
$i=i_{\lambda +n+1}.$ Thus if $\varphi$ is $\Pi_n$ and true of $G(<i),S\cap
i,A\cap i$ then $\varphi$ is forced by some coding condition  $p\in P^{<i}$
($p$ in the generic determined by $R)$ and hence by the
$(G(<\infty),S,A)-\Sigma _n$ stability of $i,$ we get that $\varphi$ is
true of $G(<\infty),S,A.$ $\dashv$

We built $R$ as in Theorem 2 by perturbing the indiscernibles. However with
extra care we can in fact obtain indiscernible preservation.

\vskip10pt

\noindent
{\bf Theorem 9.} \ There is a real $R\in L[O^\#]$ such that $R$ is generic
but not strictly generic over $L,$ $L$-cofinalities equal
$L[R]$-cofinalities and $I^R=I.$

\vskip10pt

\noindent
{\bf Proof.} \ Instead of using the $i_{\lambda +n},n\in\omega $ ($\lambda
$ limit or $0)$ use the $i^n_\alpha ,n\in\omega$ where $i^n_\alpha =$ least
element of $A_n$ greater than $i_\alpha.$ Thus $\bigcup\{i^n_\alpha
|n\in\omega \}=i_{\alpha +1}$ and as above we can construct $S$ to preserve
indiscernibles and $L$-cofinalities and satisfy that no $i^n_\alpha ,n$
odd is $S$-stable, $i^{n+1}_\alpha $ is $(S,A)-\Sigma _n$ stable
for large enough $\alpha ,n$ (given any $L$-amenable $A)$ and $A_n\cap
i_\alpha ^{n+1}$ is definable over $\lan L_i[S],S\cap i\ran$ for
$i=i_\alpha ^{n+1},$ $n$ even. Then code $(G(<\infty),S)$ by a real,
preserving indiscernibles and cofinalities, requiring as before that for
inaccessible $\kappa ,$ any coding condition with $\kappa $ in its domain
reduces dense $D\subseteq p^{<\kappa }$ strictly below $\kappa ,$ when $D$
is definable over $\lan L_\kappa [G(<\kappa ), S\cap\kappa ],$ $G(<\kappa
),S\cap\kappa \ran.$ Then for any $L$-amenable $A,$ $i^{n+1}_\alpha $ will
be $(R,A)-\Sigma _n$ stable for sufficiently large $\alpha ,n.$ This
implies as before that $R$ is not strictly generic. $\dashv$

\vskip10pt

\vskip10pt

\noindent
{\bf Remark 1.} \ A similar argument shows: \ For any $n\in\omega $ there
is a real $R\in L[O^\#]$ which is strictly generic over $L,$ yet $G$ is not
$\Sigma _n\lan L[R],R,A\ran$ whenever $R\in L[G],G$ literally generic over
$\lan L,A\ran.$ Thus there is a strict hierarchy within strict genericity,
given by the level of definability of the literally generic $G$ from the
strictly generic real.

\vskip10pt

\noindent
{\bf Remark 2.} \ The nongeneric real $R$ constructed in Friedman [??] is
strictly generic over some $L[S]$ wher $R\notin L[S].$ The same is true of
the real $R$ constructed here to satisfy Theorem 2. This leads to:

\vskip10pt

\noindent
{\bf Questions} \ (a) \ Is there a real $R\in L[O^\#],R$ not strictly
generic over any $L[S],R\notin L[S]?$  \ (b) \ Suppose $R$ is strictly
generic over $L[S],S$ generic over $L.$ Then is $R$ generic over $L?$

\vskip10pt

\begin{center}
{\bf References}\end{center}

\vskip10pt

\noindent
Beller-Jensen-Welch \ [82] \ {\it Coding the Universe,} Cambridge
University Press.

\vskip5pt

\noindent
Friedman \ [?]\  An Immune Partition of the Ordinals.

\vskip5pt

\noindent
Friedman \ [??] \ The Genericity Conjecture, JSL.

\vskip5pt

\noindent
Friedman \ [???] \   {\it Fine Structure and Class Forcing,} preliminary book
draft.

\vskip5pt

\noindent
Friedman \ [????] \ Coding without Fine Structure, to appear.

\vskip5pt

\noindent
Jech \ [?] \ {\it Set Theory,} Academic Press.

\vskip5pt

\noindent
Stanley, M. \ [?] \  Backwards Easton Forcing and $O^\#.$

\vskip5pt

\noindent
Department of Mathematics
M.I.T.
Cambridge, MA 02139
\end{document}